\documentclass[12pt,a4paper]{article}
\usepackage{amsmath,amsfonts,amssymb,amsthm,epsfig,epstopdf,titling,url,array}
\usepackage{mathtools}
\usepackage{authblk}
\usepackage[english]{babel}
\usepackage[T1]{fontenc}
\newtheorem{thm}{Theorem}[section]
\newtheorem{lem}{Lemma}[section]

\newtheorem{prop}{Proposition}[section]
\numberwithin{equation}{section}
\newcommand{\pre}{\textbf{Proof.}}
\newcommand{\R}{\mathbb{R}}

\newcommand{\C}{\mathbb{C}}

\def\fd{\hfill$\square$\par\noindent\vskip0.2cm}
\title{Bargmann transform and  generalized heat Cauchy problems}
\author{Anouar Abdelmajid Saidi, Ahmed Yahya Mahmoud and Mohamed Vall Ould Moustapha}
\begin{document}
\maketitle

\begin{abstract}
In this article we solve explicitly some Cauchy problems of the heat type attached to the generalized real and complex Dirac, Euler
 and Harmonic oscillator operators. Our principal tool is the Bargmann transform. 
\end{abstract}

\section{Introduction}
The classical Cauchy problem for the heat equation is of the form 
\begin{align}\label{HE}\qquad \left \{\begin{array}{cc}
\displaystyle  \Delta u(t, x)=\frac{\partial}{\partial t}u(t, x);\quad(t, x)\in \R^\ast_+\times \R^n,
\vspace*{0.3cm} 
\\ u(0, x)=u_0(x) ;\quad u_0 \in L^2(\R^n).
\end{array}\right.
\end{align}
where $\displaystyle\Delta=\sum_{j=1}^n\frac{\partial^2}{\partial x_j^2}$ is the classical Laplacian of $\R^n$.
We can formally replace the operator $\Delta$ by any partial differential operator and speak about generalized heat Cauchy problem. If we replace $\Delta$ by the Dirac operator $\displaystyle\frac{\partial}{\partial x}$, we obtain the famous transport equation  (see \cite{E}). 
\\ In this work, we solve explicitly some type of generalized heat Cauchy problems. Our method is based on the Bargmann transform which is an isometry from $L^2(\R)$ to the Fock space $F_\frac{a}{2}^2$ (see \cite{Z1} ). The paper is outlined as follows
\\In section 2,  old and new results concerning the Bargmann transform are given.
\\In section 3,  we compute explicitly the exact solutions of the heat Cauchy problems associated to the generalized complex and real Dirac operators 
\begin{align}\label{T1}\qquad \left \{\begin{array}{cc}
\displaystyle  D_z^a U(t, z)=\frac{\partial}{\partial t}U(t, z);\quad(t, z)\in \R^\ast_+\times \C,
\vspace*{0.3cm} 
\\ U(0, z)=U_0(z) ;\quad U_0 \in F_\frac{a}{2}^2.
\end{array}\right.\end{align}
and
\begin{align}\label{T2}\qquad \left \{\begin{array}{cc}
\displaystyle  d^a_x u(t,x)=\frac{\partial}{\partial t}u(t,x);\quad(t,x)\in \R^\ast_+\times \R,
\vspace*{0.3cm} 
\\ u(0,x)=u_0(x) ;\quad u_0 \in L^2(\R).
\end{array}\right.\end{align}
where the generalized real and complex Dirac operators are given by
\begin{align}\label{D}
  D^a_z=  \frac{1}{a}\frac{\partial }{\partial z}+\frac{z}{2}, \ \ \ \ \ \ \ \  \ \ \ \ d^a _x=  \frac{\partial }{\partial x}-ax.
\end{align}
In section 4, we solve  explicitly the following heat Cauchy problems associated to the generalized real and complex Euler operators 
\begin{align}\label{E1}\qquad \left \{\begin{array}{cc}
\displaystyle  e^a_x v(t,x)=\frac{\partial}{\partial t}v(t,x);\quad(t,x)\in \R^\ast_+\times \R,
\vspace*{0.3cm} 
\\ v(0,x)=v_0(x) ;\quad v_0 \in L^2(\R).
\end{array}\right.\end{align}
and
\begin{align}\label{E2}\qquad \left \{\begin{array}{cc}
\displaystyle E^a_z Y(t,w)=\frac{\partial}{\partial t}Y(t,w);\quad(t,w)\in \R^\ast_+\times \C,
\vspace*{0.3cm} 
\\ Y(0,w)=Y_0(w) ;\quad Y_0\in F_\frac{a}{2}^2.
\end{array}\right.
\end{align}
where the generalized real and complex Euler operators are given by
\begin{align}\label{E}
  E^a_z=  -2az\frac{\partial}{\partial z}-a, \ \ \ \ \ \ \ \  \ \ \ \ e_x^a= ax \frac{\partial }{\partial x}.
\end{align}
\noindent 
In section 5  the explicit solutions of the following heat Cauchy problem for the (real) harmonic oscillators is computed.
\begin{align}\label{H1}\qquad \left \{\begin{array}{cc}
\displaystyle  h^a_x y(t, x)=\frac{\partial}{\partial t}y(t, x);\quad(t, x)\in \R^\ast_+\times \C,
\vspace*{0.3cm} 
\\ y(0, x)=y_0(x) ;\quad y_0 \in L^2(\R).
\end{array}\right.\end{align}

\begin{align}\label{HR}
  h^a_x=\frac{\partial^2}{\partial x^2}-a^2\,x^2.
\end{align}

Section 6 is devoted to the following generalized heat Cauchy problem 
\begin{align}\label{H2}\qquad \left \{\begin{array}{cc}
\displaystyle  H^a_z V(t, z)=\frac{\partial}{\partial t}V(t, z);\quad(t, z)\in \R^\ast_+\times \C,
\vspace*{0.3cm} 
\\ V(0, z)=V_0(z) ;\quad V_0 \in F_\frac{a}{2}^2.
\end{array}\right.\end{align}
associated to the operator
\begin{align}\label{H}
  H^a_z=  \frac{\partial^2 }{\partial z^2}-\frac{a^2z^2}{4}-\frac{a}{2}, \ \ \ \ \ \ \ \  \ \ \ h^a_x=\frac{\partial^2}{\partial x^2}-a^2\,x^2.
\end{align}
called here the complex harmonic oscillator.\\

Note that the wave Cauchy and Poisson problems associated to the real Dirac, Euler and Harmonic oscillator are considered in \cite{MM1, MM2}. 
\section{Bargmann transform} 
For  $a>0$, we define a Gaussian measure on $\C$ as follows
\begin{align} d\lambda_a(z)=\frac{a}{\pi}\,e^{-a|z|^2}dz.\end{align}
The Fock space, denoted $F_a^2$, is the subspace   of all  entire functions in $L^2(\C,d\lambda_a)$ which is a Hilbert space with the inner product
 \begin{align} \left<f,g\right>=\int_\C f(z)\;\overline{g(z)}d\lambda_a(z).\end{align} 
V. Bargmann \cite{B1, B2} and G. Folland \cite{F} have defined a mapping $B$, from  the space of square integrable functions $L^2(\R)$ to the Fock space $F_\pi^2$, called Bargmann transform  
 \begin{align}\forall z\in\C,\forall f\in L^2(\R),\qquad \left[Bf\right](z)=2^{\frac{1}{4}}\,\int_{-\infty}^{\infty}f(x)\;e^{2\pi\,xz\,-\pi\,x^2-\frac{\pi\,z^2}{2}}\;dx.\end{align} 
We use in this paper a parametrized form of the Bargmann transform given by K. Zhu \cite{Z1} as follows
\begin{align}\left[B_af\right](z)=\left(\frac{2a}{\pi}\right)^{\frac{1}{4}}\,\int_{-\infty}^{\infty}f(x)\;e^{2a\,xz\,-a\,x^2-\frac{a\,z^2}{2}}\;dx.\end{align} 
This mapping is an isometry from $L^2(\R)$ to $F_a^2$, its inverse is given by 
\begin{align}\left[B_a^{-1}f\right](x)=\left(\frac{2a}{\pi}\right)^{\frac{1}{4}}\,\int_{\C} f(z)\;e^{2ax\overline{z}-ax^2-\frac{a}{2}\overline{z}^2}\,d\lambda_a(z).\end{align}

In what follows, we give some formulas involving Bargmann transform.
\begin{lem}\label{l1} (Appendix 1 of \cite{KCT}) 
For $a>0$, we have 
\begin{enumerate}
\item $\displaystyle \left[B_\frac{a}{2}(x\,f)\right](z)=\left(\frac{1}{a}\,\frac{\partial}{\partial z}+\frac{z}{2}\right)\left[B_\frac{a}{2}f\right](z)$.
\item $\displaystyle \left[B_\frac{a}{2}\left(\frac{\partial}{\partial x}\,f\right)\right](z)=\left(\frac{\partial}{\partial z}-\frac{a}{2}z\right)\left[B_\frac{a}{2}f\right](z)$.
\item $\displaystyle \left[B_\frac{a}{2}\left(\left(\frac{\partial}{\partial x}\,-ax\right)\,f\right)\right](z)=-az\,\left[B_\frac{a}{2}f\right](z)$.
\item $\displaystyle \left[B_\frac{a}{2}\left(\left(\frac{\partial}{\partial x}\,+ax\right)\,f\right)\right](z)=2\,\frac{\partial}{\partial z}\left[B_\frac{a}{2}f\right](z)$.
\end{enumerate}
\end{lem}
The following proposition is a direct consequence of lemma \ref{l1}.
\begin{prop}\label{p1}
\begin{enumerate}
\item Let $E_x^a$ and $h_x^a$ be respectively the complex Euler and real harmonic operators given in \eqref{E} and \eqref{H}. Then, we have 
$$\left[B_\frac{a}{2}(h_x^af)\right](z)=E_z^a\left[B_\frac{a}{2}f\right](z).$$
\item Let $e_x^a$ and $H_z^a$ be respectively the real Euler and complex harmonic operators given in \eqref{E} and \eqref{H}. Then, we obtain 
$$\left[B_\frac{a}{2}(e_x^af)\right](z)=H_z^a\left[B_\frac{a}{2}f\right](z).$$
\end{enumerate}
\end{prop}
The following lemma gives a new result about the Bargmann transform of the Gaussian function. 
\begin{lem}\label{th2} Let $a$, $b$, $c$ and $s$ be  real numbers such that $a>0$ and  $b>\frac{a}{2}$, then the  following lemma holds
\begin{enumerate}
\item $\displaystyle \left[B_\frac{a}{2}\left(e^{a\frac{x^2}{2}}\,e^{-b(x-s)^2}\right)\right](z)=\left(\frac{a}{\pi}\right)^\frac{1}{4}\sqrt{\frac{\pi}{b}}\,e^{asz}\,e^{\frac{a}{4}\left(\frac{a}{b}-1\right)z^2}$. 
\item   $\left[B_\frac{a}{2}\left(e^{-cx^2}\right)\right](z)=\left(\frac{a}{\pi}\right)^\frac{1}{4}\sqrt{\frac{\pi}{c+\frac{a}{2}}}\,e^{\frac{a}{4}z^2\left(\frac{a-2c}{a+2c}\right)}$. 
\end{enumerate}
\end{lem}

\noindent
\pre\,
For $s\in\R$, $a>0$ and $b>\frac{a}{2}$ we can write  
$$ \left[B_\frac{a}{2}\left(e^{a\frac{x^2}{2}}\,e^{-b(x-s)^2}\right)\right](z)=\left(\frac{a}{\pi}\right)^{\frac{1}{4}}\,\int_{-\infty}^{\infty}e^{-b(x-s)^2}\;e^{a\,xz-\frac{a\,z^2}{4}}\;dx$$
$$
=\left(\frac{a}{\pi}\right)^{\frac{1}{4}}\,\int_{-\infty}^{\infty}e^{-b(x-s-\frac{az}{2b})^2}\;e^{asz}\;e^{\frac{a^2}{4b}z^2}\,e^{-\frac{a}{4}z^2}\;dx.$$ 
From the formula (see lemma 2 of \cite{DZ})
$$\int_{-\infty}^{\infty}e^{-b(x-z)^2}\,dx=\sqrt{\frac{\pi}{b}},$$
 we obtain
$$ \left[B_\frac{a}{2}\left(e^{a\frac{x^2}{2}}\,e^{-b(x-s)^2}\right)\right](z)=\left(\frac{a}{\pi}\right)^{\frac{1}{4}}\,\sqrt{\frac{\pi}{b}}\;e^{asz}\;e^{\frac{a}{4}(\frac{a}{b}-1)z^2}\;dx.$$
 Thus we obtain the first assertion of the lemma. The second assertion is a consequence of the first one. \fd
The following lemma is well known for the classical Bargmann transform $B_\pi$ (see \cite{Z2}), we adapt its proof for the parametrized form $B_a$ in the appendix of this paper.
\\ Let  $F_r$ be a general form of the Fourier transform defined by
\begin{align}\left[F_rf\right](x)=\sqrt{\frac{ar}{\pi}}\int_{\R} f(t)\,e^{iarxt}\,dt, \quad  \textrm{ for all } a, r>0 \textrm{ and } f\in L^2(\R),\end{align} and let 
$h_r$ be the operator of $L^2(\R)$ defined as  
$$\left[h_rf\right](x)=f(r\,x),\quad \textrm{ for all }f\in L^2(\R) \textrm{ and }x\in\R.$$
\begin{lem}\label{th4} 
For all function $f\in F_\frac{a}{2}^2(\R)$, we have
\begin{enumerate}
\item $\displaystyle\left[B_\frac{a}{2}F_rB_\frac{a}{2}^{-1}f\right](z)=
2\sqrt{\frac{r}{r^2+1}}e^{-\frac{a}{4}z^2\left(\frac{r^2-1}{r^2+1}\right)}\int_{\C} f(w)e^{\frac{a}{4}\overline{w}^2(\frac{1-r^2}{r^2+1})}e^{\frac{iarz\overline{w}}{r^2+1}} \,d\lambda_\frac{a}{2}(w).$
\item $\displaystyle\left[B_\frac{a}{2}F_1B_\frac{a}{2}^{-1}f\right](z)=\sqrt{2}\,f(iz)$.
\item $\displaystyle\left[B_\frac{a}{2}F_1^{-1}B_\frac{a}{2}^{-1}f\right](z)=\frac{1}{\sqrt{2}}\,f(-iz)$.
\item $\displaystyle\left[B_\frac{a}{2}\,h_r\,B_\frac{a}{2}^{-1}f\right](x)=\sqrt{\frac{2}{r^2+1}}e^{-\frac{a}{4}z^2\left(\frac{r^2-1}{r^2+1}\right)}\int_{\C} f(-iw)e^{\frac{a}{4}\overline{w}^2(\frac{1-r^2}{r^2+1})}e^{\frac{iarz\overline{w}}{r^2+1}} \,d\lambda_\frac{a}{2}(w).$
\end{enumerate}
\end{lem}
For recent works on Bargman transform see \cite{CHH, DZ}.
\section{Heat Cauchy problems for the complex and real Dirac operators}
\begin{thm}\label{th6}
The solution of  the heat Cauchy problem \eqref{T1} associated to the complex Dirac operator is given by the formula 
\begin{align}U(t,z)=e^\frac{zt}{2}\,e^{\frac{t^2}{4a}}\,U_0\left(z+\frac{t}{a}\right).\end{align}
 \end{thm}
\noindent
\pre\,
Let $U$  be a solution of \eqref{T1}.
Then, by applying the inverse of the Bargmann transform $B^{-1}_\frac{a}{2}$, we obtain
 $$\qquad \left \{\begin{array}{cc}
\displaystyle  B^{-1}_\frac{a}{2}\left(\left(\frac{1}{a} \frac{\partial }{\partial z}+\frac{z}{2}\right)U(t,z)\right)=\frac{\partial}{\partial t}B^{-1}_\frac{a}{2}(U(t,z));\quad(t, z)\in \R^\ast_+\times \C,
\vspace*{0.3cm} 
\\B^{-1}_\frac{a}{2}(U(0, z))=B^{-1}_\frac{a}{2}U_0;\quad U_0 \in F_\frac{a}{2}^2.
\end{array}\right.$$
Using lemma \ref{l1} we get 
$$\qquad \left \{\begin{array}{cc}
\displaystyle x\,\left[B^{-1}_\frac{a}{2}U\right](t,x)=\frac{\partial}{\partial t}\left[B^{-1}_\frac{a}{2}U\right](t,x);\quad(t,x)\in \R^\ast_+\times\R,
\vspace*{0.3cm} 
\\\left[B^{-1}_\frac{a}{2}U\right](0,x)=\left[B^{-1}_\frac{a}{2}U_0\right](x) ;\quad U_0\in F_\frac{a}{2}^2.
\end{array}\right.$$ 
 Thus $B^{-1}_\frac{a}{2}U$ satisfies the formula
 $$\left[B^{-1}_\frac{a}{2}U\right](t,x)=e^{tx}\,\left[B^{-1}_\frac{a}{2}U_0\right](x).$$
 This implies that 
 $$U(t,z)=\left[B_\frac{a}{2}\left(e^{tx}B^{-1}_\frac{a}{2}U_0\right)\right](z)$$
$$=\left(\frac{a}{\pi}\right)^{\frac{1}{2}}\,\int_{\R}e^{a\,xz\,-\frac{a}{2}\,x^2-\frac{a\,z^2}{4}}\;e^{tx}\int_{\C} \,U_0(v)\; e^{ax\overline{v}-\frac{a}{2}x^2-\frac{a}{4}\overline{v}^2}\,d\lambda_\frac{a}{2}(v)\,dx$$
$$=\left(\frac{a}{\pi}\right)^{\frac{1}{2}}\,\int_{\R}\int_{\C} e^{tx}\,U_0(v)\;e^{a\,xz\,-\frac{a}{2}\,x^2-\frac{a\,z^2}{4}} e^{ax\overline{v}-\frac{a}{2}x^2-\frac{a}{4}\overline{v}^2}\,d\lambda_\frac{a}{2}(v)\,dx$$
$$=\left(\frac{a}{\pi}\right)^{\frac{1}{2}}\,e^{-\frac{a\,z^2}{4}} e^{\frac{a\,(z+\frac{t}{a})^2}{4}}\int_{\R}\int_{\C} U_0(v)\;e^{a\,x(z+\frac{t}{a})\,-\frac{a}{2}\,x^2}e^{-\frac{a\,(z+\frac{t}{a})^2}{4}} e^{ax\overline{v}-\frac{a}{2}x^2-\frac{a}{4}\overline{v}^2}\,d\lambda_\frac{a}{2}(v)\,dx$$
$$=e^\frac{zt}{2}\,e^{\frac{t^2}{4a}}\left[B_\frac{a}{2}(B^{-1}_\frac{a}{2}U_0)\right]\left(z+\frac{t}{a}\right)=e^\frac{zt}{2}\,e^{\frac{t^2}{4a}}\,U_0\left(z+\frac{t}{a}\right).$$\fd
\begin{thm}\label{th7}
The solution of  the heat Cauchy problem \eqref{T2} associated to the real Dirac operator is given by the formula
\begin{align}u(t,x)=e^{-axt}\,e^{\frac{-at^2}{2}}\,u_0(x+t).\end{align}
 \end{thm}

\noindent
\pre\,
Let $u$  be a solution of \eqref{T2}.
Then, by applying  the Bargmann transform $B_\frac{a}{2}$, we obtain
 $$\left \{\begin{array}{cc}
\displaystyle  B_\frac{a}{2}\left(\left(\frac{\partial }{\partial x}-ax\right)u(t,x)\right)=\frac{\partial}{\partial t}B_\frac{a}{2}\left(u(t, x)\right);\quad(t,x)\in \R^\ast_+\times \R, 
\vspace*{0.3cm} 
\\ B_\frac{a}{2}(u(0,x))=B_\frac{a}{2}(u_0) ;\quad u_0 \in L^2(\R).
\end{array}\right.$$
Using lemma \ref{l1} we obtain  
$$\qquad \left \{\begin{array}{cc}
\displaystyle -az\,\left[B_\frac{a}{2}u\right](t,z)=\frac{\partial}{\partial t}\left[B_\frac{a}{2}u\right](t,z);\quad(t,z)\in \R^\ast_+\times\C,
\vspace*{0.3cm} 
\\\left[B_\frac{a}{2}u\right](0,z)=\left[B_\frac{a}{2}u_0\right](z) ;\quad u_0\in L^2(\R).
\end{array}\right.$$ 
 Then $B_\frac{a}{2}u$ satisfies the formula 
 $$\left[B_\frac{a}{2}u\right](t,z)=e^{-azt}\,\left[B_\frac{a}{2}u_0\right](z).$$
 This implies that 
 $$u(t,x)=\left[B^{-1}_\frac{a}{2}(e^{-azt}\,B_\frac{a}{2}u_0)\right](x)$$
$$=\left(\frac{a}{\pi}\right)^{\frac{1}{2}}\,\int_{\C}\int_{\R} e^{-azt}u_0(s)\;e^{a\,sz\,-\frac{a}{2}\,s^2-\frac{a\,z^2}{4}} e^{ax\overline{z}-\frac{a}{2}x^2-\frac{a}{4}\overline{z}^2}\,ds\;d\lambda_a(z)$$
$$=\left(\frac{a}{\pi}\right)^{\frac{1}{2}}\,\int_{\C}\int_{\R} u_0(s)\;e^{az(s-t)\,-\frac{a}{2}\,s^2-\frac{a\,z^2}{4}} e^{ax\overline{z}-\frac{a}{2}x^2-\frac{a}{4}\overline{z}^2}\,ds\;d\lambda_a(z)$$
By setting  $s'=s-t$, we obtain 
 $$u(t,x)=\left(\frac{a}{\pi}\right)^{\frac{1}{2}}\,\int_{\C}\int_{\R} u_0(s'+t)\;e^{azs'}e^{-\frac{a}{2}\,s'^2}e^{-\frac{a}{2}\,t^2}e^{-a\,s't}e^{-\frac{a\,z^2}{4}} e^{ax\overline{z}-\frac{a}{2}x^2-\frac{a}{4}\overline{z}^2}\,ds'\,d\lambda_a(z)$$
 $$=\left(\frac{a}{\pi}\right)^{\frac{1}{2}}\,\int_{\C} e^{ax\overline{z}-\frac{a}{2}x^2-\frac{a}{4}\overline{z}^2}\int_{\R} \bigg(u_0(s'+t)\,e^{-\frac{a}{2}\,t^2}e^{-a\,s't}\bigg)\;e^{azs'}e^{-\frac{a}{2}\,s'^2}e^{-\frac{a\,z^2}{4}} \,ds'\,d\lambda_a(z)$$
$$=\left[B^{-1}_\frac{a}{2}\bigg(B_\frac{a}{2}\bigg(u_0(s'+t)\,e^{-\frac{a}{2}\,t^2}e^{-a\,s't}\bigg)\bigg)\right](x)=u_0(x+t)\,e^{-\frac{a}{2}\,t^2}\,e^{-a\,xt}.$$\fd
\section{Heat Cauchy problems for the complex and  real Euler operators }
\begin{thm}\label{th8}
The solution of  the heat Cauchy problem \eqref{E1} associated to the real Euler operator is given by the formula
 \begin{align}v(t,x)=v_0(e^{at}\,x).\end{align}
\end{thm}
\noindent
\pre\,
Let $e$ be the solution of \eqref{E1}. We define the function $y(t,x)$ by
$$y(t,x)=e(t,e^{-at}x).$$ We get easily  $\frac{\partial}{\partial t}y(t,x)=0$ then 
$y(t,x)=y(0,x)=e(0,x)=e_0(x).$ Then 
$$e(t,e^{-at}x)=e_0(x),\textrm{ for all } (t,x)\in \R_+\times \R.$$
Consequently, $e(t,x)=e_0(e^{at}\,x)$ for all $(t,x)\in \R_+\times \R.$\fd
\begin{thm}\label{th9}
The solution of  the heat Cauchy problem \eqref{E2} associated to the complex Euler operator is given by the formula
 \begin{align}Y(t,z)=e^{-at}\,Y_0(e^{-2at}z).\end{align}
\end{thm}

\noindent
\pre\,
The proof of the theorem is Mutatis mutandis  the proof of theorem \ref{th8}.\fd
\section{Heat Cauchy problems for the real harmonic oscillators}
In this section, we study the heat Cauchy problem for the real harmonic oscillator $\displaystyle h^a_x=\frac{\partial^2}{\partial x^2}-a^2\,x^2.$ It is clear that  $$h_x^a=\left(\frac{\partial}{\partial x}\,-ax\right)\left(\frac{\partial}{\partial x}\,+ax\right).$$

Note that the heat kernel for the classical real harmonic oscillator  has been known for a long time as Mehler kernel  (\cite{BGV} p.145) and is given by the formula
$$K_a(x,s,t)=$$
\begin{align}\label{r00}\sqrt{\frac{a}{\pi}}\frac{1}{\sqrt{\sinh(2at)}}\;\,{\rm exp}\left\{-\frac{a}{2}(x^2+s^2)\coth(2at)+\frac{axs}{\sinh(2at)}\right\},\end{align}
but the method used here and the obtained formula are new. 
\begin{thm}\label{th10}
The heat Cauchy problem \eqref{H1}   associated to the real harmonic oscillator has the unique solution given by 
$$y(t,x)=\int_{\R} K_a(x,s,t)\;y_0(s)\;ds$$
Where $K_a(x,s,t)=$
\begin{align}\label{r0}\sqrt{\frac{a}{\pi}}\,(e^{2at}-e^{-2at})^{-\frac{1}{2}}\,{\rm exp}\left\{-\frac{a(e^{at}x-e^{-at}s)^2}{e^{2at}-e^{-2at}}+\frac{a}{2}(x^2-s^2)\right\}.\end{align}
\end{thm}

\noindent
\pre\, Let $\gamma$  be a solution of \eqref{H1}. Then, by applying the Bargmann transform $B_\frac{a}{2}$, we obtain
 $$\qquad \left \{\begin{array}{cc}
\displaystyle \left[B_\frac{a}{2}(h_x^a\gamma(t,x))\right](z)=\frac{\partial}{\partial t}\left[B_\frac{a}{2}(\gamma(t,x))\right](z);\quad(t,z)\in \R^\ast_+\times \C,
\vspace*{0.3cm} 
\\\left[B_\frac{a}{2}(\gamma(0,x))\right](z)=\left[B_\frac{a}{2}\gamma_0\right](z) ;\quad \gamma_0\in L^2(\R).
\end{array}\right.$$
Proposition \ref{l1} assures that 
$$\qquad \left \{\begin{array}{cc}
\displaystyle \left(-2az\frac{\partial}{\partial z}-a\right)\left[B_\frac{a}{2}(\gamma(t,x))\right](z)=\frac{\partial}{\partial t}\left[B_\frac{a}{2}(\gamma(t,x))\right](z);\quad(t,z)\in \R^\ast_+\times\C,
\vspace*{0.3cm} 
\\\left[B_\frac{a}{2}(\gamma(0,x))\right](z)=\left[B_\frac{a}{2}(\gamma_0)\right](z) ;\quad \gamma_0\in L^2(\R).\end{array}\right.$$
Now, theorem \ref{th9} gives 
$$\left[B_\frac{a}{2}\gamma\right](t,z)=e^{-at}\left[B_\frac{a}{2}\gamma_0\right](e^{-2at}z),\qquad\forall (t,z)\in\R_+^\ast\times\C.$$
So $$\gamma(t,x)=\left[B_\frac{a}{2}^{-1}\left(e^{-at}\left[B_\frac{a}{2}\gamma_0\right](e^{-2at}z)\right)\right](x)$$
$$=\sqrt{\frac{\alpha}{\pi}}e^{-at}\int_{\C} e^{ax\overline{z}-\frac{a}{2}x^2-\frac{a}{4}\overline{z}^2}\int_{\R} \gamma_0(s)e^{ase^{-2at}z-\frac{a}{2}s^2-\frac{a}{4}e^{-4at}z^2}\,ds\,d\lambda_\frac{a}{2}(z).$$
Using the change of variables $s'=e^{-2at}s$, we get 
$$\gamma(t,x)=\sqrt{\frac{a}{\pi}}\;e^{at}\int_{\C} e^{ax\overline{z}-\frac{a}{2}x^2-\frac{a}{4}\overline{z}^2}\int_{\R} \gamma_0(e^{2at}s')e^{as'z-\frac{a}{2}e^{4at}s'^2-\frac{a}{4}e^{-4at}z^2}\,ds'\,d\lambda_\frac{a}{2}(z)$$
$$=\sqrt{\frac{a}{\pi}}\,e^{at}\int_{\R} \gamma_0(e^{2at}s')\,e^{-\frac{a}{2}e^{4at}s'^2} \int_{\C} e^{as'z-\frac{a}{4}e^{-4at}z^2+ax\overline{z}-\frac{a}{2}x^2-\frac{a}{4}\overline{z}^2}\,d\lambda_\frac{a}{2}(z)\,ds'$$
$$=\left(\frac{a}{\pi}\right)^\frac{1}{4}\,e^{at}\int_{\R} \gamma_0(e^{2at}s')\,e^{-\frac{a}{2}e^{4at}s'^2}\left[B_\frac{a}{2}^{-1}\left(e^{as'z-\frac{a}{4}e^{-4at}z^2}\right)\right](x)\,ds'.$$
We deduce from lemma \ref{th2} (with $\frac{a}{b}-1=-e^{-4at}$) that 
$$\left[B_\frac{a}{2}^{-1}\left(  e^{as'z-\frac{a}{4}e^{-4at}z^2}\right)\right](x)=\left(\frac{\pi}{a}\right)^\frac{1}{4}\frac{\sqrt{a}}{\sqrt{\pi}\sqrt{1-e^{-4at}}}e^{\frac{a}{2}x^2}\,e^{\frac{-a(x-s')^2}{1-e^{-4at}}}.$$ So we obtain 
$$\gamma(t,x)=e^{at}e^{\frac{a}{2}x^2}\frac{\sqrt{a}}{\sqrt{\pi}\sqrt{1-e^{-4at}}}\int_{\R} \gamma_0(s'e^{2at})\,e^{-\frac{a}{2}e^{4at}s'^2}\,e^{\frac{-a(x-s')^2}{1-e^{-4at}}}ds'.$$
Set again $s=e^{2at}s'$, we get 
$$\gamma(t,x)=e^{-at}e^{\frac{a}{2}x^2}\frac{\sqrt{a}}{\sqrt{\pi}\sqrt{1-e^{-4at}}}\int_{\R} \gamma_0(s)\,e^{-\frac{a}{2}s^2}\,e^{\frac{-a(x-se^{-2at})^2}{1-e^{-4at}}}ds$$
$$=\sqrt{\frac{a}{\pi}}\,\frac{1}{\sqrt{e^{2at}-e^{-2at}}}\;\int_{\R} \gamma_0(s)\,e^{\frac{a}{2}(x^2-s^2)}\,e^{\frac{-a\left(xe^{at}-se^{-at}\right)^2}{e^{2at}-e^{-2at}}}ds$$
$$=\int_{\R} K_a(x, s, t)\gamma_0(s)ds.$$ \fd
The formula \eqref{r0} obtained in the  theorem \ref{th10} agree with the well known formula given in \eqref{r00}. In fact,
$$K_a(x,s, t)=\sqrt{\frac{a}{\pi}}\,(e^{2at}-e^{-2at})^{-\frac{1}{2}}\,{\rm exp}\left\{\frac{-a(e^{at}x-e^{-at}s)^2}{e^{2at}-e^{-2at}}+\frac{a}{2}(x^2-s^2)\right\}$$
$$=\sqrt{\frac{a}{2\pi}}\frac{1}{\sqrt{\sinh(2at)}}\,{\rm exp}\left\{\frac{-a(e^{2at}x^2+e^{-2at}s^2-2xs)}{e^{2at}-e^{-2at}}+\frac{a}{2}x^2-\frac{a}{2}s^2\right\}$$
$$=\sqrt{\frac{a}{2\pi}}\frac{1}{\sqrt{\sinh(2at)}}\,{\rm exp}\left\{-\frac{a}{2}x^2\left(\frac{2e^{2at}}{e^{2at}-e^{-2at}}-1\right)\,-\frac{a}{2}s^2\left(\frac{2e^{-2at}}{e^{2at}-e^{-2at}}+1\right)\,+\frac{2axs}{e^{2at}-e^{-2at}}\right\}$$
$$=\sqrt{\frac{a}{2\pi}}\frac{1}{\sqrt{\sinh(2at)}}\,{\rm exp}\left\{-\frac{a}{2}(x^2+s^2)\left(\frac{e^{2at}+e^{-2at}}{e^{2at}-e^{-2at}}\right)+\frac{axs}{\sinh(2at)}\right\}.$$
Thus we obtain the formula \eqref{r00}.\fd
\section{Heat Cauchy problems for the complex harmonic oscillators}
In this section, we study the heat Cauchy problem for the complex harmonic oscillator $\displaystyle H^a_z=\frac{\partial^2 }{\partial z^2}-\frac{a^2z^2}{4}-\frac{a}{2}$. It is clear that  $$H^a_z=\left(\frac{\partial}{\partial z}\,+az\right)\left(\frac{\partial}{\partial z}\,-az\right).$$
\begin{thm}\label{th11}
The Cauchy problem \eqref{H2}  for the heat equation associated to the complex harmonic oscillator has the unique solution given by 
$$V(t,z)=\int_{\C} J_a(z,w',t)\,V_0(w')\;\,d\lambda_\frac{a}{2}(w'),$$
where \begin{align} J_a(z,w',t)=\frac{2i}{\sqrt{\cosh(at)}}\;exp\left\{\frac{a}{4}(\overline{w'}^2-z^2)\tanh(at) +\frac{az\overline{w'}}{2\cosh(at)}\right\}. \end{align}
 \end{thm}
\noindent
\pre\,
Let $E$ be the solution of \eqref{H2}. By applying the inverse of the Bargmann transform $B^{-1}_\frac{a}{2}$, we obtain 
 $$\qquad \left \{\begin{array}{cc}
\displaystyle B^{-1}_\frac{a}{2}\displaystyle \left(\left(\frac{\partial }{\partial z}+\frac{az}{2}\right)\left(\frac{\partial }{\partial z}-\frac{az}{2}\right)E(t, z)\right)=\frac{\partial}{\partial t}B^{-1}_\frac{a}{2}(E(t,z));\quad(t, z)\in \R^\ast_+\times \C,
\vspace*{0.3cm} 
\\ B^{-1}_\frac{a}{2}(E(0, z))=B^{-1}_\frac{a}{2}E_0 ;\quad E_0 \in F_\frac{a}{2}^2.
\end{array}\right.$$
Using Proposition \ref{l1} we obtain  
$$\qquad \left \{\begin{array}{cc}
\displaystyle ax\,\frac{\partial}{\partial x}\left[B^{-1}_\frac{a}{2}E\right](t,x)=\frac{\partial}{\partial t}\left[B^{-1}_\frac{a}{2}E\right](t,x);\quad(t,x)\in \R^\ast_+\times\R,
\vspace*{0.3cm} 
\\\left[B^{-1}_\frac{a}{2}E\right](0,x)=\left[B^{-1}_\frac{a}{2}E_0\right](x) ;\quad U_0\in F_\frac{a}{2}^2.
\end{array}\right.$$ 
So  theorem \ref{th8} assures that 
 $$\left[B^{-1}_\frac{a}{2}E\right](t,x)=\left[B^{-1}_\frac{a}{2}E_0\right](e^{at}\,x).$$
 This implies that 
 $$E(t,z)=\left[B_\frac{a}{2}\bigg(\left[B^{-1}_\frac{a}{2}E_0\right](e^{at}\,x)\bigg)\right](z).$$
By Lemma \ref{th4} we obtain
$$E(t,z)=\sqrt{\frac{2}{r^2+1}}e^{-\frac{a}{4}z^2\left(\frac{r^2-1}{r^2+1}\right)}\int_{\C} f(-iw)e^{\frac{a}{4}\overline{w}^2(\frac{1-r^2}{r^2+1})}e^{\frac{iarz\overline{w}}{r^2+1}} \,d\lambda_\frac{a}{2}(w),\textrm{ where } r=e^{at}.$$
Then $$E(t,z)=\frac{2}{\sqrt{\cosh(at)}}\;e^{-\frac{a}{4}z^2\;\tanh(at)}\int_{\C} f(-iw)\;e^{-\frac{a}{4}\overline{w}^2\tanh(at)}\;e^{\frac{iaz\overline{w}}{2\cosh(at)}}\,d\lambda_\frac{a}{2}(w).$$
Using the change of variable $w'=-iw$ we get 
$$E(t,z)=\;\int_{\C} f(w')\;\frac{2i}{\sqrt{\cosh(at)}}\;e^{\frac{a}{4}(\overline{w'}^2-z^2)\tanh(at)} \;e^{\frac{az\overline{w'}}{2\cosh(at)}}\,d\lambda_\frac{a}{2}(w').$$\fd
\newpage\pagestyle{plain}
\section*{Appendix :} 
\addcontentsline{toc}{section}{Appendix :}
We give here  the proof of lemma \ref{th4}. 

\noindent 
\pre\,\textbf{of lemma \ref{th4}}.
\begin{enumerate}
\item Firstly, we have
$$
F_r\left(B_\frac{a}{2}^{-1}f\right)(w)=\left(\frac{a}{\pi}\right)^{\frac{1}{4}}\sqrt{\frac{ar}{\pi}}\int_{\R} e^{iarxt}\int_{\C} f(w) e^{at\overline{w}-\frac{a}{2}t^2-\frac{a}{4}\overline{w}^2}\,d\lambda_\frac{a}{2}(w)\,dt$$
$$=\left(\frac{a}{\pi}\right)^{\frac{1}{4}}\sqrt{\frac{ar}{\pi}}\int_{\C} f(w)e^{-\frac{a}{4}\overline{w}^2}\int_{\R}  e^{\frac{-a}{2}(t-irx-\overline{w})^2}\,e^{\frac{a}{2}(irx+\overline{w})^2}dt\,d\lambda_\frac{a}{2}(w)$$
$$=\frac{\sqrt{2ar}}{(a\pi)^{\frac{1}{4}}}\int_{\C} f(w)e^{-\frac{a}{4}\overline{w}^2}\,e^{\frac{a}{2}(irx+\overline{w})^2}d\lambda_\frac{a}{2}(w)$$
$$=\frac{\sqrt{2ar}}{(a\pi)^{\frac{1}{4}}}e^{\frac{-ar^2x^2}{2}}\int_{\C} f(w)e^{\frac{a}{4}\overline{w}^2}\,e^{airx\overline{w}}d\lambda_\frac{a}{2}(w).$$
Thus we obtain that \,$(B_\frac{a}{2}F_rB_\frac{a}{2}^{-1}f)(z)=$
$$\left(\frac{a}{\pi}\right)^{\frac{1}{4}}\frac{\sqrt{2ar}}{(a\pi)^{\frac{1}{4}}}\int_{\R} e^{axz-\frac{a}{2}x^2-\frac{a}{4}z^2} e^{\frac{-ar^2x^2}{2}}\int_{\C} f(w)e^{\frac{a}{4}\overline{w}^2}\,e^{airx\overline{w}}d\lambda_\frac{a}{2}(w)\,dx$$
$$=\sqrt{\frac{2ar}{\pi}}e^{-\frac{a}{4}z^2}\int_{\C} f(w)e^{\frac{a}{4}\overline{w}^2}\,\int_{\R} e^{-x^2(\frac{a}{2}+\frac{ar^2}{2})}e^{x(az+iar\overline{w})} \,dx\,d\lambda_\frac{a}{2}(w).$$
We denote $c=\frac{a}{2}+\frac{ar^2}{2}$. Using the change of variable $t=\sqrt{c}.\,x$ we get
$$(B_\frac{a}{2}F_rB_\frac{a}{2}^{-1}f)(z)$$
$$=\sqrt{\frac{2ar}{c\,\pi}}e^{-\frac{a}{4}z^2}\int_{\C} f(w)e^{\frac{a}{4}\overline{w}^2}\,\int_{\R} e^{-t^2+\frac{t}{\sqrt{c}}\left(az+iar\overline{w}\right)} \,dt\,d\lambda_\frac{a}{2}(w)$$
$$=\sqrt{\frac{2ar}{c}}e^{-\frac{a}{4}z^2}\int_{\C} f(w)e^{\frac{a}{4}\overline{w}^2}\,\int_{\R} e^{-\left(t-\frac{az+iar\overline{w}}{2\sqrt{c}}\right)^2}e^{\left(\frac{az+iar\overline{w}}{2\sqrt{c}}\right)^2} \,dt\,d\lambda_\frac{a}{2}(w)$$
$$=\sqrt{\frac{2ar}{c\,\pi}}e^{-\frac{a}{4}z^2}\int_{\C} f(w)e^{\frac{a}{4}\overline{w}^2}e^{\left(\frac{az+iar\overline{w}}{2\sqrt{c}}\right)^2} \,d\lambda_\frac{a}{2}(w)$$
$$=\sqrt{\frac{2ar}{c}}e^{-\frac{a}{4}z^2\left(1-\frac{a}{c}\right)}\int_{\C} f(w)e^{\frac{a}{4}\overline{w}^2(1-\frac{ar^2}{c})}e^{\frac{ia^2zr\overline{w}}{2c}} \,d\lambda_\frac{a}{2}(w)$$
$$=2\sqrt{\frac{r}{r^2+1}}e^{-\frac{a}{4}z^2\left(\frac{r^2-1}{r^2+1}\right)}\int_{\C} f(w)e^{\frac{a}{4}\overline{w}^2(\frac{1-r^2}{r^2+1})}e^{\frac{iarz\overline{w}}{r^2+1}} \,d\lambda_\frac{a}{2}(w).$$
\item In the case  $r=1$, we obtain 
 $$(B_\frac{a}{2}F_1B_\frac{a}{2}^{-1}f)(z)=\sqrt{2}\int_{\C} f(w)e^{\frac{iaz\overline{w}}{2}} \,d\lambda_\frac{a}{2}(w).$$ By the formula of reproducing kernel (see  proposition 2.2 of \cite{Z1}) this function is exactly equal to  $\sqrt{2}\,f(iz)$. 
\item On the other hand, $$(B_\frac{a}{2}F_1^{-1}B_\frac{a}{2}^{-1})=(B_\frac{a}{2}F_1B_\frac{a}{2}^{-1})^{-1},$$ which justifies the assertion 3.
\item  A simple calculus gives $\sqrt{r}\,h_r\circ F_1=F_r \textrm{ and thus } h_r=\frac{1}{\sqrt{r}}\,F_r\circ F_1^{-1}$, then
$$(B_\frac{a}{2}\,h_r\,B_\frac{a}{2}^{-1}f)(z)=\frac{1}{\sqrt{r}}\,(B_\frac{a}{2}\,F_r\,F_1^{-1}\,B_\frac{a}{2}^{-1}f)(z)=\frac{1}{\sqrt{r}}\,(B_\frac{a}{2}\,F_r\,B_\frac{a}{2}^{-1}\,B_\frac{a}{2}\,F_1^{-1}\,B_\frac{a}{2}^{-1}f)(z).$$
 Using assertion 1 and 3, we deduce that
$$(B_\frac{a}{2}\,h_r\,B_\frac{a}{2}^{-1}f)(z)=\sqrt{\frac{2}{r^2+1}}e^{-\frac{a}{4}z^2\left(\frac{r^2-1}{r^2+1}\right)}\int_{\C} f(-iw)e^{\frac{a}{4}\overline{w}^2(\frac{1-r^2}{r^2+1})}e^{\frac{iarz\overline{w}}{r^2+1}} \,d\lambda_\frac{a}{2}(w).$$
\end{enumerate}
\fd

\newpage 
\newcommand{\Addresses}{{
  \bigskip
  \footnotesize

 Anouar Saidi, \textsc{Department of Mathematics,
 College of Arts and Sciences-Gurayat,
 Jouf University-Kingdom of Saudi Arabia.}\par\nopagebreak
 \textsc{Department of Mathematics, Faculty of Sciences of Monastir,5019 Monastir-Tunisia}.
\\\textit{E-mail address}: \texttt{saidi.anouar@yahoo.fr}
  \medskip
 
 A. Y. Mahmoud,\textsc{Department of Mathematics,
 College of Arts and Sciences-Gurayat,
 Jouf University-Kingdom of Saudi Arabia.}\par\nopagebreak
 \textsc{Shendi University ,Faculte of science and Technologe, Departement of Mathematics ,Shendi Sudan}.
  \\\textit{E-mail address}: \texttt{humudi999@gmail.com}

  \medskip

  M.V. Ould Moustapha, \textsc{Department of Mathematic,
 College of Arts and Sciences-Gurayat,
 Jouf University-Kingdom of Saudi Arabia }\par\nopagebreak
\textsc{ Faculte des Sciences et Techniques
Universit\'e de  Nouakchott Al-Aasriya.
Nouakchott-Mauritanie}\\
  \textit{E-mail address}, M. V.~Ould Moustapha: \texttt{mohamedvall.ouldmoustapha230@gmail.com}

}}

\Addresses

\end{document}